\documentclass[12pt]{amsart}

\usepackage{amsmath}
\usepackage{latexsym}
\usepackage{amsfonts}
\usepackage{amssymb}
\newtheorem{thm}{Theorem}[section]
\newtheorem{lem}[thm]{Lemma}
\newtheorem{prop}[thm]{Proposition}

\theoremstyle{definition}
\newtheorem*{defi}{Definition}
\newtheorem*{obs}{Observation}
\newtheorem*{note}{Note}

\newcommand{\R}{\mathbb R}
\newcommand{\C}{\mathbb C}
\newcommand{\Z}{\mathbb Z}
\newcommand{\N}{\mathbb N}

\begin{document}

\title{Model spaces results for the Gabor and Wavelet transforms.}%

\author{Gerard Ascensi}
\thanks{The author is supported by MTM2005-08984-C02-01 and 2005SGR00611 projects and by the Fundaci\'{o} Cr\`{e}dit Andorr\`{a}.}
\address{Faculty of Mathematics, University of Vienna, Nordbergstrasse 15, 1090 Vienna, Austria}
\email{gerard.ascensi@univie.ac.at}


\begin{abstract}
We prove that the unique Gabor atom with analytical model space is
the gaussian function. We give an analogous result for the wavelet
transform. For the general case we give a new approach to study
the irregular Gabor and wavelet frames. We improve some results
for Gabor atoms in the Feichtinger algebra, and for a special
class of wavelets.
\end{abstract}

\maketitle

\section{Introduction}

Gabor and wavelet analysis are two topics in harmonic analysis
that have had and important growth in the last 30 years. This is
due to easy ways of constructing frames that these two transforms
give and their applications in  several fields. Most of the
examples of Gabor and wavelet frames correspond to consider
regular nets of points. That is, in the Gabor case sets of the
type $\{e^{-2\pi ibnt}g(t-am)\}_{n,m\in \Z}$ and in the wavelet
case of type
$\{2^{-\frac{j}{2}}\psi(2^{-j}t-n)\}_{j,n\in\Z}$. There
are plenty of references about this topic (see e.g. \cite{Dau} or
\cite{Mal}). We can find sufficient and necessary conditions for
the existence of this kind of frames, and also a big number of
applications, both theoretical as practical. At the same time
interest arose in the study of the frames that come from an
irregular grid, study that can give a better knowledge of the
structure and properties of these transforms.

Thus, the problem that we want to study here is when a set of functions
$\{e^{-2\pi b_nt}g(t-a_n)\}$ with $\Lambda=\{(a_n,b_n)\}$ a
discrete set of the plane can be a frame of $L^2(\R)$. Also, the
equivalent problem can be proposed for the wavelet case. We will
give the appropriate definitions later. Classical harmonic
analysis does not deal with this kind of problems due to the
irregularity of the grid. Instead, functional analysis and Hilbert
spaces techniques must be used.

First approaches can be found in \cite{FG1} and \cite{FG2}. By
means of an abstract theory of representations these papers give a
way to study Gabor and wavelet irregular frames. They prove, in
case  the window function  belongs to the so called Feichtinger
algebra, that if every small ball of the plane contains a point of
$\Lambda$, then it will define a frame. Feichtinger and Grochenig
used functional analysis, and for this reason they do not obtain
specific bounds. Regarding more sufficient conditions, in last
years  a series of papers has appeared (\cite{SZ1}, \cite{SZ2},
\cite{FS1}, \cite{FS2} \cite{SZ3} and \cite{SZ4}), giving the same kind of
results,  with explicit bounds, also under special conditions on
the functions.

The aim of this paper is to give another approach to this problem,
based in the methods of \cite{OS1}. Fixed an analyzing function,
we will consider the model space of all continuous transforms, so
that our problem becomes a sampling problem in the space. These
sets (sampling and also the dual notion of interpolation) have
been studied in the case of spaces of holomorphic functions, as
for example Hardy or Paley-Wiener spaces. For some analyzing
functions the model space is a  space of holomorphic functions.
The most important example is the Gaussian function, for which the
model space can be identified with the Fock space, in which the
sampling and interpolation seta are completely characterized (see
\cite{Lyu}, \cite{S4}, and \cite{SW1}). Unfortunately, a first
result of this paper states that this is essentially the only
possible example. For the wavelet case we achieve a similar
result, in this case we prove that, under some restrictions, the
only analyzing wavelet for which the model space consists of
holomorphic functions is the Poisson wavelet. In this case the
spaces are the Bergman spaces, for which also exists a
characterization of the sampling sets \cite{S2}. We will give the
detailed statements and proofs in section 3.

In section 4 we obtain sufficient conditions analogous to the
previously mentioned ones, for the Gabor transform, using the
ideas in \cite{OS1} introduced in the wavelet case. We restrict to
analyzing windows in the Feichtinger algebra, and prove first
that any set that is closed to a sampling set is also sampling,
and that a sufficiently dense set is a sampling set. Both results
are known, but the proof we give leads to explicit bounds. With
this method we are able to extend some results that up to now were
only known in case the model space consists of holomorphic
functions. In section 5 we review briefly the wavelet situation
and fix a small gap in \cite{OS1}.

We do not deal here with necessary conditions that the discrete
set must satisfy in order to define a frame (see \cite{RS1} and
\cite{CDH1} for the Gabor case and \cite{HK1} and \cite{HK2}
for the wavelet case). We refer to \cite{H1} and \cite{K1}
for a general overview.

\section{Definitions and preliminary results.}

Fixed a function $g\in L^2(\R)$ with $\|g\|=1$, recall that the
windowed Fourier transform (also called Gabor transform) of $f\in
L^2(\R)$ is defined through correlation with the atoms
$g_z(t)=e^{-2\pi iyt}g(t-x), z=x+iy$,
\begin{equation*}
Gf(z)=\langle f(t),g_z(t)\rangle=\int_{\R}f(t)\overline{e^{-2\pi
iyt}g(t-x)}\,dt.
\end{equation*}

One can prove that the reconstruction formula
\begin{equation}\label{formularecongabor}
f(t)=\int_{\R}\int_{\R}Gf(z)g_z(t)\,dm(z)\quad \forall f\in
L^2(\R)
\end{equation}
holds in the $L^2$-sense, and in particular,
\begin{equation*}
\|f\|^2=\langle f,f\rangle=\int_{\R}\int_{\R}Gf(z)\langle g_z,f\rangle
\,dm(z)=\|Gf\|^2.
\end{equation*}

Similarly, given $\psi\in L^2(\R), \|\psi\|=1$ the {\bf wavelet
transform} of $f$ is defined in the upper half-plane, $z=x+iy,
y>0$ using correlation with the wavelets
$\psi_z(t)=y^{\frac{-1}{2}}\psi\left(\frac{t-x}{y}\right)$,
\begin{equation*}
Wf(z)=\langle f,\psi_z\rangle=\int_{-\infty}^{\infty}f(t)
\overline{y^{\frac{-1}{2}}\psi\left(\frac{t-x}{y}\right)}\,dt.
\end{equation*}

As before, one has a continuous reconstruction formula, this time
when $\psi$ is {\bf admissible \bf} in the sense that
\begin{equation*}
\int_0^{\infty}\frac{|\widehat{\psi}(\xi)|^2}{\xi}\,d\xi=1.
\end{equation*}

The reconstruction formula  is
\begin{equation}\label{formularecon}
f(t)=\int_0^{\infty}\int_{-\infty}^{\infty}Wf(z)\psi_z(t)\,
d\mu(z)
\end{equation}
with $d\mu(z)=\frac{dx\,dy}{y^2}$, to be interpreted again in the
$L^2$-sense, in particular
\begin{equation*}
\|f\|^2= \int_0^{\infty}\int_{-\infty}^{\infty}|Wf(z)|^2 \,
d\mu(z)=\|Wf\|_{\R\times\R^+}^2.
\end{equation*}

Next we recall that  a family  $\{f_k\}_{k\in\N}$ in a Hilbert
space $H$ is said to be a {\bf frame} if there are constants $A,B>
0 $ such that
\begin{equation*}
A\|f\|^2\leqslant\sum_{k=1}^{\infty}|\langle f,f_k\rangle|^2
\leqslant B\|f\|^2,\quad \forall f\in H.
\end{equation*}
The constants $A$ and $B$ are called the frame bounds. In this
case, there is a dual frame $\{\widetilde{f}_k\}_{k\in\N}$ with
frame bounds $\frac{1}{B}$ and $\frac{1}{A}$ such that
\begin{equation*}
f=\sum_{k=1}^{\infty}\langle f,\widetilde{f}_k\rangle f_k \quad
\\
f=\sum_{k=1}^{\infty}\langle f,f_k\rangle \widetilde{f}_k \quad
\forall f\in H
\end{equation*}
(see details in \cite{Dau}).

\begin{defi}
Let $\Lambda=\{z_n\}_{n\in\N}$ be a discrete set in $\C$ and $g\in
L^2(\R)$ a window function with $\|g\|=1$.  We define the {\bf
Gabor system} $G(g,\Lambda)$ as:
\begin{equation*}
G(g,\Lambda)=\bigl\{g_z(t)=e^{-2\pi iyt}g(t-x);
z=x+iy\in\Lambda\bigr\}.
\end{equation*}

Similarly, if $\Sigma=\{\sigma_n\}_{n\in\N}$  is a discrete set in
$\R\times\R^+$ and $\psi\in L^2(\R)$ is admissible function we
define the {\bf wavelet system} $W(\psi,\Sigma)$ as:
\begin{equation*}
W(\psi,\Sigma)=\left\{\psi_{\sigma}(t)=y^{\frac{-1}{2}}\psi\left(\frac{t-x}{y}\right);
\sigma=x+iy\in\Sigma\right\}.
\end{equation*}
\end{defi}

We are interested in knowing when these systems are frames of
$L^2(\R)$. A restatement is in terms of a  sampling property in
the space of all of transforms,
\begin{equation*}
H_g=\left\{F\in L^2(\C) : \exists f\in L^2(\R): F(z)=Gf(z)=\langle
f,g_z\rangle\right\}
\end{equation*}
\begin{footnotesize}
\begin{equation*}
H_{\psi}=\left\{F\in L^2(\R\times\R^+): \exists f\in L^2(\R):
F(z)=Wf(z)=\langle f,\psi_z\rangle\right\}.
\end{equation*}
\end{footnotesize}
We call them \emph{model spaces}, and their description is as
follows (see \cite{Dau}). The model space $H_g$  is a Hilbert
subspace of $L^2(\C)$ that is characterized by the  reproductive
kernel:
\begin{equation*}
k_g(z,z_0)=k(z,z_0)=k_{z_0}(z)=\langle g_{z_0},g_z\rangle.
\end{equation*}
That is, $F\in H_g$ if and only if $F\in L^2(\C)$ and
\begin{equation}\label{reproducciogabor}
F(z_0)=\int_{\C}F(z)\overline{k(z,z_0)}\,dm(z).
\end{equation}

Similarly, the model space $H_{\psi}$ of an admissible wavelet is
a Hilbert subspace of $L^2(\R\times\R^+, d\mu)$  characterized by
the  reproductive kernel
\begin{equation*}
k_{\psi}(z,z_0)=k(z,z_0)=k_{z_0}(z)=\langle
\psi_{z_0},\psi_z\rangle.
\end{equation*}
That is, $F\in H_{\psi}$ if and only if $F\in L^2(\R\times\R^+,
d\mu)$ and
\begin{equation}\label{reproduccioondetes}
F(z_0)=\int_{\R\times\R^+}F(z)\overline{k(z,z_0)}\,d\mu(z).
\end{equation}

In terms of $k_g(z)=\langle g,g_z\rangle=k_g(0,z)$ one has
\begin{align*}
k_g(z,z_0)=&\overline{\langle
g_z,g_{z_0}\rangle}=\overline{e^{-2\pi ix(y-y_0)}\langle
g,g_{z_0-z}\rangle}\\
=&\overline{e^{-2\pi ix(y-y_0)}k_g(z_0-z)}
\end{align*}
and hence the reproduction formula \eqref{reproducciogabor} takes
the form:
\begin{equation*}
F(z_0)=\int_{\C}F(z)e^{-2\pi ix(y-y_0)}k_g(z_0-z)\,dx\,dy.
\end{equation*}
This is called a \emph{twisted convolution} and, correspondingly,
$F_{z_0}(z)=e^{2\pi ix_0(y-y_0)}F(z-z_0)$ is called a
\emph{twisted translation}. Since twisted translation can be seen
to be a continuous operation in $L^2(\C)$, it follows easily that
all $F\in H_g$ are uniformly continuous; more precisely, given
$\varepsilon$ there exists $\delta$ such that if
$|z_1-z_2|<\delta$, one has
\begin{equation*}
\bigl||F(z_1)|-|F(z_2)|\bigr|<\|F\|\varepsilon=\|f\|\varepsilon.
\end{equation*}

In the wavelet setting, the equation \eqref{reproduccioondetes} is
also a convolution, but with respect to the hyperbolic (or affine)
group. If we define $z_0\cdot z=y_0z+x_0$ (where the juxtaposition
is the usual product of $\C$) we can endow a group structure to
$\R\times\R^+$. In this group the identity will be $i=(0,1)$ and
$z_0^{-1}\cdot z=(z-z_0)/y_0$. This group is not commutative. Its
left-invariant measure (for translations with respect to the
group) is precisely $d\mu(z)=\frac{dx\,dy}{y^2}$. In terms of
$k_{\psi}(z)=\langle \psi,\psi_z\rangle=k_{\psi}(0,z)$,
\begin{equation*}
k_{\psi}(z,z_0)=\langle \psi_{z_0},\psi_z\rangle=\langle
\psi_{z^{-1}\cdot z_0},\psi\rangle=\overline{k_{\psi}(z^{-1}\cdot
z_0)}
\end{equation*}
and so  \eqref{reproduccioondetes} is
\begin{equation*}
F(z_0)=\int_{\R\times\R^+}F(z)k(z^{-1}\cdot z_0)\,d\mu(z)
\end{equation*}
where we find a formula of non commutative convolution. Again,
these functions are uniformly continuous; in this case is more
natural to use the hyperbolic distance in the half-plane
\begin{equation*}
d(z_1,z_2)=\frac{1}{2}\log\frac{1+\overline{d}(z_1,z_2)}{1+\overline{d}(z_1,z_2)}
\end{equation*}
where $\overline{d}(z_1,z_2)$ is the pseudohyperbolic distance:
\begin{equation*}
\overline{d}(z_1,z_2)=\left|\frac{z_1-z_2}{z_1-\overline{z_2}}\right|.
\end{equation*}

The precise statement is that given $\varepsilon>0$ there exists
$\delta>0$ such that if $d(z_1,z_2)<\delta$, for every $F\in H$,
one has
$|F(z_1)-F(z_2)|\leqslant\|F\|\varepsilon=\|f\|\varepsilon$.

\begin{defi}
A discrete set $\Lambda=\{\lambda_n\}_{n\in\Z}\subset\C$ is said
to be a {\bf sampling set} for $H_g$ if there are constants
$A,B>0$ such that
\begin{equation*}
A\|F\|^2\leqslant\sum_{n\in\Z}|F(\lambda_n)|^2\leqslant
B\|F\|^2\quad \forall F\in H_g.
\end{equation*}

Similarly, a discrete set
$\Sigma=\{\sigma_n\}_{n\in\Z}\subset\R\times\R^+$ is said to be a
{\bf sampling set} for $H_{\psi}$ if there are constants $A,B>0$
such that
\begin{equation*}
A\|F\|^2\leqslant\sum_{n\in\Z}|F(\sigma_n)|^2\leqslant
B\|F\|^2\quad \forall F\in H_{\psi}.
\end{equation*}
\end{defi}

With these definitions, it is clear that   $G(g,\Lambda)$ is a
frame of $L^2(\R)$ if and only if $\Lambda$ is a sampling set for
$H_g$ and $W(\psi,\Sigma)$ is a frame of $L^2(\R)$ if and only if
$\Sigma$ is a sampling set for $H_{\psi}$.

We have to restrict the choice of the analyzing function. This is
the same restriction that is used in \cite{FG1} and \cite{FG2}.
Here we find one of the most important differences between the
Gabor and the wavelet case.

\begin{defi}
We define the {\bf Feichtinger algebra} as the set of functions
$g\in L^2(\R)$ such that
\begin{equation*}
k(z)=\langle g,g_z\rangle=\int_{\R}g(t)e^{-2\pi
ity}\overline{g(t-x)}\,dt\in L^1(\C).
\end{equation*}
We will denote this set as $\mathcal{A}$.
\end{defi}

\begin{defi}
Given a continuous function $F$ defined in $\C$ we define its {\bf
local maximal function} as:
\begin{equation*}
MF(z)=\sup_{|w-z|<1} |F(w)|.
\end{equation*}
\end{defi}

The windows that one needs to consider are those for which the
maximal function of the reproductive kernel is integrable.
Surprisingly this happens if $g\in\mathcal{A}$ (see \cite{FG2}

\begin{prop}
Let $g\in L^2(\R)$ be a Gabor atom and $k$ the reproductive kernel
of its model space. If $g\in\mathcal{A}$ then $Mk$ is integrable
in $\C$. That is, $Mk$ is integrable if $k$ is.
\end{prop}

In the wavelet case we need considering an analogous notion:

\begin{defi}
We define the set of {\bf wavelets with integrable kernel} as
those admissible wavelets $\psi\in L^2(\R)$ such that
\begin{equation*}
k(z)=\langle
\psi,\psi_z\rangle=\int_{\R}\psi(t)y^{-\frac{1}{2}}\overline{\psi\bigl(\frac{t-x}{y}\bigr)}\,dt\in
L^1(\R\times\R^+).
\end{equation*}
We denote the set of all the wavelets with integrable kernel as
$\mathcal{B}$.
\end{defi}

This set preserves some of the properties of the Feichtinger
algebra, but not all of them. This is because in the Gabor case
the group that plays a role is the Weyl-Heisenberg group, that is
unimodular, while in the wavelet case is the affine group, that is
not. As an example, $Gf(z)$ is integrable if and only if $Gf(-z)$
is but in the wavelet setting $Wf(z)$ can be integrable and
$Wf(z^{-1})$ not. This is one of the differences of the structure
of the group.

\begin{defi}
Given a continuous function $F$ defined in $\R\times\R^+$ we
define its {\bf local maximal function} as:
\begin{equation*}
MF(z)=\sup_{w\in B(z,1)} |F(w)|
\end{equation*}
where $B(z,1)$ denotes the ball of center $z$ and radius $1$ in
$\R\times\R^+$ using the hyperbolic distance. The area of this
ball for the left invariant measure is $4\pi
\sinh^2(\frac{1}{2})$.
\end{defi}

The set of wavelets such that the local maximal function of its
kernel is integrable differs from $\mathcal{B}$ (see
\cite{Gro}), and we will call it $\mathcal{MB}$.

\section{Analytic model spaces.}

A good knowledge of the model spaces is very convenient for the
study of this kind of frames. The best situation occurs when the
model space is a space of holomorphic functions, in which case
there are some complete characterizations. We begin with the Gabor
case.

\begin{defi}
The {\bf Fock space} $\mathcal{F}$ is defined:
\begin{equation*}
\left\{ F \text{ entire with }
\|F\|_{\mathcal{F}}^2=\int_{\C}|F(z)|^2e^{-\pi|z|^2}\,dm(z)<\infty\right\}.
\end{equation*}
\end{defi}

The Bargmann transform gives the relationship between this space
and $L^2(\R)$.

\begin{defi}
Given $f\in L^2(\R)$, its {\bf Bargmann transform} of $f$ is

\begin{equation*}
Bf(z)=2^{\frac{1}{4}}\int_{\R}f(t)e^{2\pi tz-\pi
t^2-\frac{\pi}{2}z^2}\,dt.
\end{equation*}
\end{defi}

This transform is an isomorphism between $L^2(\R)$ and
$\mathcal{F}$ \cite{Fol}. The relation between the Bargmann and
the Gabor transform using the gaussian function is given by

\begin{equation*}
Gf(x+iy)=e^{-\frac{\pi}{2}|x+iy|^2}e^{\pi ixy}Bf(x+iy).
\end{equation*}

In an informal way we can think that the model space of the
Gaussian function is the Fock space; in fact, it is easily checked
that they have the same sampling and interpolation sets. In the
Fock space we have a complete characterization of these sets. This
characterization has been proved in an independent way in
\cite{Lyu}, \cite{S4} and \cite{SW1}.

The case of the Gaussian function is exceptional because it is the
only atom with this property:

\begin{thm} Consider the model space of a Gabor atom $g\in L^2(\R)$
\begin{equation*}
H_g=\{F(z)=\int_{\R}f(t)e^{2\pi ity}g(t-x)\,dt,\; f\in L^2(\R)\}.
\end{equation*}
Then this space is a space of holomorphic functions, modulo a
multiplication by a weight, if and only if $g$ is a time-frequency
translation of the Gaussian function.
\end{thm}

\begin{proof}
We suppose that there is a $M(z)=M(x,y)$ such that
\begin{equation*}
M(z)F(z)\in Hol(\C)\qquad \forall F\in\overline{H}.
\end{equation*}
Then
\begin{align*}
\overline{\partial}MF(z)=&\int_{\R}f(t)\overline{\partial}\bigl(M(z)e^{2\pi
ity}g(t-x)\bigr)\,dt\\
=&0\qquad \forall f\in L^2(\R)
\end{align*}
which is equivalent to
\begin{equation*}
\overline{\partial}\bigl(M(z)e^{2\pi
ity}g(t-x)\bigr)=0\qquad\forall t\in\R.
\end{equation*}
The left hand side equals
\begin{align*}
(\overline{\partial}M) e^{2\pi ity}g(t-x)-&\frac 12 M(z) e^{2\pi
ity}g'(t-x)\\
-&\pi t M(z)e^{2\pi ity}g(t-x)
\end{align*}
so we obtain
\begin{equation*}
2\bigl(\overline{\partial}M(x,y)-M(x,y)\pi
t\bigl)g(t-x)=M(x,y)g'(t-x).
\end{equation*}

Changing to variable $w=t-x$ we get the differential equation
\begin{equation*}
\frac{g'(w)}{g(w)}=2\frac{\overline{\partial}M}{M}-2\pi (w+x).
\end{equation*}

This forces $\frac{\overline{\partial}M}{M}-\pi x$ to be a
constant $c_1$ and implies $\frac{g'(w)}{g(w)}=2 c_1-2\pi w$,
which in turn gives that $g$ is a time-frequency translation of
the Gaussian function.
\end{proof}

\begin{note}
This result is also true if we assume that the model space
consists of quasiregular functions, that is, satisfying an
equation of type $\overline{\partial} (MF)=\mu \partial{MF}$ with
some dilation factor $\mu$.
\end{note}

To obtain examples of wavelets with analytical model space one has
to change a bit the definition of the wavelet transform. We have
to consider this transform defined in the Hardy space $H^2(\R)$
instead of $L^2(\R)$. There is not substantial difference between
both cases.  For $\alpha>1$, we define the {\bf Poisson wavelet}
$\psi^{\alpha}(t)$ as:
\begin{equation*}
\psi^{\alpha}(t)=c_{\alpha}(t+i)^{-\frac{\alpha+1}{2}}
\end{equation*}
where $c_{\alpha}$ is a normalizing constant. The model space of
these wavelets are in one to one correspondence with the Bergman
spaces.

\begin{defi}
For $\alpha>1$, we define the {\bf Bergman space} of the
half--plane $A_{\alpha}(\R\times\R^+)$ as:
\begin{align*}
\Bigl\{F \text{ analytical in
}&\R\times\R^+:\\
&\quad |F\|_{\alpha}^2=\int_{\R^\times\R^+}|F(z)|^2y^{\alpha}\,d\mu(z)<\infty\Bigr\}.
\end{align*}
\end{defi}

It is easy to see that, modulo a  multiplication by
$y^{\frac{\alpha}{2}}$,  both spaces have the same reproducing
kernel. For the Bergman spaces there exists as well a
characterization of the sampling sets. We can find it in \cite{S2}
for the equivalent versions on the disk.

Now we give a uniqueness theorem for the wavelet case. We will
prove more generally,  for harmonic function spaces. But we have
to put a restriction in the weight function.

\begin{thm}
Let $\psi$ be an admissible and real valued wavelet. Then there is
a weight $\omega(y)$ such that each function $\omega(y)Wf(x,y)$
for $f\in L^2(\R)$ is harmonic if and only if $\psi$ is a linear
combination of $\Re (t+i)^{\alpha}$ and $\Im (t+i)^{\alpha}$ with
$\alpha<-1$.
\end{thm}

\begin{proof}
We have to study when $\Delta\omega(y)Wf(x,y)=0$:
\begin{equation*}
\Delta\omega(y)Wf(x,y)=\int_{-\infty}^{\infty}f(t)\Delta
y^{\frac{-1}{2}}\omega(y)\psi\Bigl(\frac{t-x}{y}\Bigr)\,
dt=0
\end{equation*}
for every $f\in L^2(\R)$. This is equivalent to:
\begin{equation*}
\Delta y^{\frac{-1}{2}}\omega(y)\psi\Bigl(\frac{x}{y}\Bigr)=0.
\end{equation*}

Denoting $\beta(y)=y^{\frac{-1}{2}}\omega(y)$ we obtain
\begin{multline*}
\Delta\left(\beta(y)\psi\Bigl(\frac{x}{y}\Bigr)\right)=2\frac{x}{y^3}\beta(y)\psi'\Bigl(\frac{x}{y}\Bigr)
+\beta''(y)\psi\Bigl(\frac{x}{y}\Bigr)\\
-2\frac{x}{y^2}\beta'(y)\psi'\Bigl(\frac{x}{y}\Bigr)+\Bigl(\frac{x^2}{y^2}+1\Bigr)\frac{1}{y^2}
\beta(y)\psi''\Bigl(\frac{x}{y}\Bigr)=0.
\end{multline*}

After a change of variable and multiplying the former equation by
$y^2$ we arrive to a differential equation
\begin{multline}\label{har2}
2\beta(y)t\psi'(t)+ y^2\beta''(y)\psi(t)-2y\beta'(y)t\psi'(t)\\
+\beta(y)(t^2+1)\psi''(t)=0.
\end{multline}

We introduce the change $y=e^u$ (we remember that $y>0$) and
$\gamma(u)=\beta(e^u)$. Then
$\gamma'(u)=e^u\beta'(e^u)=y\beta'(y)$,
$\gamma''(u)=e^u\beta'(e^u)+\left(e^u\right)^2\beta''(e^u)=y\beta'(y)+y^2\beta''(y)$
and we can write \eqref{har2} as
\begin{multline*}
\gamma'(u)\psi(t)-\gamma'(u)\left(\psi(t)+2t\psi'(t) \right)
+\gamma(u)\\
\left((t^2+1)\psi''(t)+2t\psi'(t)\right)=0.
\end{multline*}
We can think this last equation as a scalar product in $\R^3$. If
we define $a(u)=\bigl(\gamma''(u),\gamma'(u),\gamma(u)\bigr)$,
\begin{equation*}
b(t)=\left(\psi(t),-\psi(t)-2t\psi'(t),(t^2+1)\psi''(t)+2t\psi'(t)\right)
\end{equation*}
we can write it
\begin{equation}\label{harpe}
\langle a(u),b(t)\rangle=0.
\end{equation}

We call $E$ to the subspace of $\R^3$ generated by $a(u)$ when $u$
varies, and $F$ to the one generated by $b(t)$ when $t$ varies.
Equation \eqref{harpe} say that $E\bot F$. As both spaces are
contained in $\R^3$ we can say that $\dim E+\dim F\leqslant 3$.
This inequality gives a very small number of options. First we
have the trivial cases $\dim F=0$, that corresponds to
$\psi(t)=0$, and $\dim E=0$ that corresponds to
$\gamma(u)=0\Rightarrow\omega(y)=0$.

If $\dim F=1$, this means that there exist a vector
$(c_1,c_2,c_3)$ and a function $B(t)$ such that $
b(t)=B(t)(c_1,c_2,c_3)$, that is,
\begin{equation*}
\left.\begin{aligned} \psi(t)&=c_1B(t)\\
-\psi(t)-2t\psi'(t)&=c_2B(t)\\
(t^2+1)\psi''(t)+2t\psi'(t)&=c_3B(t)
\end{aligned}\right\}.
\end{equation*}
This implies
\begin{equation*}
\frac{\psi(t)}{-\psi(t)-2t\psi'(t)}=\frac{c_1}{c_2}.
\end{equation*}

The cases where the denominator cancels are trivial. Solving the
differential equation we obtain $ \psi(t)=\mathrm{C}e^{\alpha
t^2}$, but this will not be an admissible wavelet for real values
of $\alpha$.

If $\dim E=1$, there exist a vector $(c_1,c_2,c_3)$ and a function
$A(u)$ such that $ a(u)=A(u)(c_1,c_2,c_3)$ that is,
\begin{gather*}
\left.\begin{aligned} \gamma''(u)&=A(u)c_1\\
\gamma'(u)&=A(u)c_2\\ \gamma(u)&=A(u)c_3
\end{aligned}\right\}\Rightarrow
\frac{\gamma(u)}{\gamma'(u)}=\frac{c_3}{c_2}.
\end{gather*}
This means that $\log
\gamma(u)=\left(e^u\right)^{\frac{c_3}{c_2}}\mathrm{C}$. This is
equivalent to $\beta(y)=y^{\frac{c_3}{c_2}}=y^{\alpha}$. For this
$\beta$, using \eqref{har2} we see that $\psi$ satisfies
\begin{equation*}
2t\psi'(t)+\alpha(\alpha-1)\psi(t) -2\alpha
t\psi'(t)+\left(t^2+1\right)\psi''(t)=0.
\end{equation*}

It is easy to check that $\psi(t)=(t+i)^{\alpha}$ satisfies this
equation.  However we search real solutions. Then we have to take
$\Re (t+i)^{\alpha}$ and $\Im (t+i)^{\alpha}$, that, by the
linearity will also be solutions of the equation. As this equation
has order $2$ and we have two solutions, we can say that the rest
of solutions are a linear combination of those. As we are only
interested in admissible wavelets of $L^2(\R)$, we have to
restrict to the case $\alpha<-1$.
\end{proof}

\section{Sampling results for the Gabor transform.}

In this section we will give some sufficient conditions for
irregular Gabor frames and obtain some stability properties. As
mentioned in the introduction, we use the techniques of
\cite{OS1}; similar results are obtained in \cite{FS1} and
\cite{FS2}. Here and in the next section we make a unified
presentation, valid both for the Gabor and the wavelet transforms,
and improve some of the results.

We will use that $Mk\in L^1(\C)$, a fact that in the Gabor case,
as pointed out before, is a consequence of $g\in\mathcal{A}$.

We establish first some notation. $H$ will always be the model
space of a normalized Gabor atom ($\|g\|=1$), and
$\Gamma=\{z_j\}_{j\in\N}$ a discrete set in $\C$.

\begin{defi}
 We will say that
$\Gamma$ is an \emph{uniformly discrete set} if
there is $\varepsilon>0$ such that $|z_i-z_j|>\varepsilon \;
\forall i\neq j$. The supremum of such $\varepsilon$ is called the
\emph{separation constant} of $\Gamma$.
\end{defi}

\begin{prop}\label{gunioseparades} If $\Gamma$ is a sampling set
for $H$ then $\Gamma$ is a finite union of uniformily discrete sets.
\end{prop}

\begin{proof}
As $\Gamma$ is a sampling set, $\exists C>0$ such that
\begin{equation}\label{gsamplingsep}
\sum_{\Gamma}|F(z_j)|^2\leqslant C\|F\|^2.
\end{equation}

If $\Gamma$ is not a finite union of uniformly discrete sets then
for every $N$ and every $\delta$ there exists a ball $B$ of radius
$\delta$ such that $|\Gamma\cap B|>N$. As $|k_{z_0}(z)|$ is an
uniformly continuous function, there exists $\delta$ such that
\begin{equation*}
\bigl||k_{z_0}(z_1)|-|k_{z_0}(z_2)|\bigr| <1/2,
\end{equation*}
if $|z_1-z_2|<\delta$, where $\delta$ does not depend on $z_0$. We
take this $\delta$ and $N>4C$. Let $B$ be the corresponding ball
of radius $\delta$ that contains $N$ points of the set. Let $w$ be
the center of $B$ and apply \eqref{gsamplingsep} to the function
$k_w(z)\in H$:
\begin{equation*}
\sum_{\Gamma}|k_w(z_j)|^2\leqslant C\|k_w\|^2=C.
\end{equation*}
But on the other hand
\begin{equation*}
\sum_{\Gamma}|k_w(z_j)|^2 \geqslant\sum_{\Gamma\cap B}
(1-1/2)^2\geqslant N\frac{1}{4}>C
\end{equation*}
and we have a contradiction.
\end{proof}

\begin{lem}\label{gacotdiscretKenL1}
Let $k$ be such that $Mk\in L^1(\C)$ and $\Lambda\subset\C$ a
uniformly discrete set with separation constant $\varepsilon$. Then
\begin{equation*}
\sum_{\lambda\in\Lambda}|k(\lambda)|<\frac{\varepsilon^{-2}}{4\pi}\|Mk\|_1
\end{equation*}
\end{lem}

\begin{proof}
We suppose without losing generality that $\frac{\varepsilon}{2}<1$. Then
\begin{equation*}
\sum_{\lambda\in\Lambda}|k(\lambda)|\leqslant\sum_{\lambda\in\Lambda}\frac{1}{|B(\lambda,\frac{\varepsilon}{2})|}
\int_{B(\lambda,\frac{\varepsilon}{2})}Mk(z)\,dm(z).
\end{equation*}

Since by hypothesis those balls are disjoint the lemma follows.
\end{proof}
\begin{obs}
Notice that the bound only depends on the separation constant of $\Lambda$.
\end{obs}

\begin{prop}\label{gcotasuperior}
Given an analyzing Gabor atom $g\in\mathcal{A}$ and $\Gamma$ a
uniformily discret set, there exists $B>0$ such that
\begin{equation*}
\sum_{\gamma\in\Gamma}|F(\gamma)|^2\leqslant B\|F\|^2\quad \forall
F\in H.
\end{equation*}
\end{prop}

\begin{proof}
Calculating directly we have that:
\begin{small}
\begin{align*}
\sum_{\gamma\in\Gamma}&|F(\gamma)|^2=\sum_{\gamma\in\Gamma}\left
|\int_{\C}F(z)\overline{k_{\gamma}(z)}\,dm(z)\right|^2\\
\leqslant&
\sum_{\gamma\in\Gamma}\left(\int_{\C}|F(z)|^2|k_{\gamma}(z)|
\,dm(z) \right)\left(\int_{\C}|k_{\gamma}(z)|\,dm(z)\right)\\
=&\int_{\C}|F(z)|^2\sum_{\gamma\in\Gamma}|k(z-\gamma)|\,dm(z)
\int_{\C}|k(z)|\,dm(z)\leqslant B\|F\|^2.
\end{align*}
\end{small}
Here we have used that $\int_{\C}|k(z)|\,dm(z)<\infty$ because the
kernel is integrable and also
\begin{equation*}
\sum_{\gamma\in\Gamma}|k(z-\gamma)|=\sum_{\lambda\in(z-\Gamma)}|k(\lambda)|
\end{equation*}
is bounded independently of $z$, since and $z-\Gamma$ has the same
separation constant than $\Gamma$ and we can apply
\ref{gacotdiscretKenL1}.
\end{proof}

\begin{thm}\label{gestab_sampling_teo}
Let $H$ be de model space of a Gabor atom $g\in\mathcal{A}$.
Given $\Lambda=\{z_j\}_{j\in\N}$ a sampling set for $H$ there
exists $\delta>0$ such that if $\Gamma=\{w_j\}_{j\in\N}$ satisfies
$|z_j-w_j|<\delta\;\forall j$, then $\Gamma$ is also a
sampling set.
\end{thm}

This result is deduced in a trivial way from the following pair of
lemmas.

\begin{lem}\label{gestab_sampling1}
Let $\Lambda=\{z_j\}_{j\in\N}$ and $\Gamma=\{w_j\}_{j\in\N}$ be
two discrete sets in $\C$, then
\begin{equation*}
\Biggl|\biggl(\sum_{j\in\N}|F(z_j)|^2\biggr)^{\frac{1}{2}}-
\biggl(\sum_{j\in\N}|F(w_j)|^2\biggr)^{\frac{1}{2}}\Biggr| \leqslant
d_1d_2\|F\|\, \forall F\in H,
\end{equation*}
where $d_1$ and $d_2$ are defined by:
\begin{itemize}
\item $d_1^2=\sup_j\int_{\C}|k_{z_j-w_j}(z)-k(z)|\,dm(z)$

\item $d_2^2=\sup_z\sum_{j\in\N}|k_{z_j}(z)-k_{w_j}(z)|$
\end{itemize}
\end{lem}

\begin{proof}
Calculating directly we have that:
\begin{multline*}
\Biggl|\biggl(\sum_{j\in\N}|F(z_j)|^2\biggr)^{\frac{1}{2}}-\biggl(\sum_{j\in\N}|F(w_j)|^2\biggr)^{\frac{1}{2}}\Biggr|\\
=\bigl|\|(F(z_j))_j\|_2-\|(F(w_j))_j\|_2\bigr|,
\end{multline*}
where $\|\cdot\|_2$ means the $l^2$--norm, thinking
$\bigl(F(z_j)\bigr)_j$ as a succession. By the triangle inequality and the reproduction formula,
\begin{align*}
\bigl|\|(F(z_j)&)_j\|_2-\|(F(w_j))_j\|_2\bigr|\\
\leqslant&
\bigl\|(|F_{-w_j}(z_j-w_j)|-|F_{-w_j}(0)|)_j\bigr\|_2\\
\leqslant&\sum_{j\in\N}\left(\int_{\C}|F_{-w_j}(z)||k_{z_j-w_j}(z)-k(z)|\,dm(z)\right)^2.
\end{align*}

Now we use the Schwartz inequality to bound the above by:
\begin{multline*}
\sum_{j\in\N}\int_{\C}|F_{-w_j}(z)|^2
|k_{z_j-w_j}(z)-k(z)|\,dm(z)\\
\cdot\int_{\C}|k_{z_j-w_j}(z)-k(z)|\,dm(z)
\end{multline*}
and use the definitions of $d_1,d_2$.
\end{proof}

\begin{lem}\label{gdspetites}
Let $\Lambda=\{z_j\}_{j\in\N}$ be a uniformly discrete set, and assume
that $g\in\mathcal{A}$. Then, for every $\varepsilon$ there exists
$\delta$ such that if $\Gamma=\{w_j\}_{j\in\N}$ satisfies
$|z_j-w_j|\leqslant\delta\;\forall j$ then $d_1d_2\leqslant
\varepsilon$, with $d_1$ and $d_2$ defined as in
\ref{gestab_sampling1}.
\end{lem}

\begin{proof}
First we will see that we can make $d_1$ as small as we want if
$\Gamma$ is close enough to $\Lambda$. To see this we write:
\begin{multline*}
\int_{\C}|k_{z_j-w_j}(z)-k(z)|\,dm(z)\\
=\int_{\C}\Bigl|e^{2\pi
i(x_j-a_j)(y-y_j+b_j)}k\bigl(z-(z_j-w_j)\bigr)-k(z)\Bigr|\,dm(z).
\end{multline*}

As $k$ is integrable, there is $R>0$ such that, if $|\alpha|<1$,
\begin{equation*}
\int_{\C\setminus B(0,R)}|k(z-\alpha)|dm(z)<\frac{\varepsilon}{4}.
\end{equation*}

In this way, if we assume $\delta<1$, we have that:
\begin{footnotesize}
\begin{align*}
\int_{\C}&\Bigl|e^{2\pi i(x_j-a_j)(y-y_j+b_j)}k\bigl(z-(z_j-w_j)\bigr)-k(z)\Bigr|\,dm(z)\\
=&\int_{B(0,R)}\Bigl|e^{2\pi i(x_j-a_j)(y-y_j+b_j)}k\bigl(z-(z_j-w_j)\bigr)-k(z)\Bigr|\,dm(z)\\
+&\int_{\C\setminus B(0,R)}\Bigl|e^{2\pi i(x_j-a_j)(y-y_j+b_j)}k\bigl(z-(z_j-w_j)\bigr)-k(z)\Bigr|\,dm(z)\\
\leqslant& \int_{B(0,R)}\Bigl|e^{2\pi i(x_j-a_j)(y-y_j+b_j)}-1\Bigr|\Bigl|k\bigl(z-(z_j-w_j)\bigr)\Bigr|\,dm(z)\\
+&\int_{B(0,R)}\Bigl|k\bigl(z-(z_j-w_j)\bigr)-k(z)\Bigr|\,dm(z)+\frac{\varepsilon}{2}.
\end{align*}
\end{footnotesize}

For the first integral, we use that $|x_j-a_j|,|y_j-b_j|\leqslant
|z_j-w_j|$. As $|y|<R$, choosing $\delta$ small enough we can
achieve $\bigl|e^{2\pi
i(x_j-a_j)(y-y_j+b_j)}-1\bigr|<\frac{\varepsilon}{4\|k\|_1}$ for
every $|y|<R$. For the second integral it is only necessary to use
that the translation operator is a continuous in $L^1(\C)$.

We go now to bound $d_2$. If $\alpha=\sup_{i\neq j} |z_i-z_j|$ is
the separation constant of $\Lambda$, we assume
$\delta<\frac{\alpha}{3}$, so that
$|w_i-w_j|>\frac{\alpha}{3}\;\forall i\neq j$. That is,
$\frac{\alpha}{3}$ can be used as separation constant for
$\Lambda$ as well as for $\Gamma$. As a matter of fact,
$\Lambda_z=\{z-z_j \}_{j\in\N}$ and $\Gamma_z=\{z-w_j\}_{j\in\N}$
also have the same separation constant. Here we can apply
\ref{gacotdiscretKenL1} to prove that there is
$C=C(\frac{\alpha}{3})$ such that $d_2\leqslant 2C$.

Hence we see that $d_2$ is bounded by $C$ when we take $\delta$
small, and $d_1$ can be made as small as we wish taking $\delta$
small enough. Therefore there exists $\delta$ such that if
$|z_j-w_j|\leqslant \delta\;\forall j$, then $d_1d_2<\varepsilon$.
\end{proof}

In proposition \ref{gunioseparades} we saw that every sampling set
is a finite union of uniformly discrete sets; now the previous lemmas allow
to improve this result.

\begin{thm}\label{gsamplingseparada}
Let $\Gamma=\{z_j\}_{j\in\N}$ be a sampling set for the model
space $H$ of a Gabor atom $g\in\mathcal{A}$. Then $\Gamma$
contains a subset $\widetilde{\Gamma}\subseteq\Gamma$ such that
$\widetilde{\Gamma}$ is a sampling and uniformly discrete set.
\end{thm}

\begin{proof}
As $\Gamma$ is sampling set we know by \ref{gunioseparades} that
it is a finite union of uniformly discrete sets, and that there are $A,B>0$
such that for every $F\in H$,
\begin{equation*}
A\|F\|^2\leqslant \sum_{j\in\N}|F(z_j)|^2\leqslant B\|F\|^2.
\end{equation*}

For every $\delta$ (we think in small $\delta$) we can
define $\widetilde{\Gamma}\subseteq\Gamma$ so that if
$\widetilde{\Gamma}=\{w_i\}_{i\in\N}$, one has that
\begin{itemize}
\item $B(w_i,\delta)\cap\widetilde{\Gamma}=\{w_i\}$

\item $\cup_{i\in\N}\bigl(B(w_i,\delta)\cap\Gamma\bigr)=\Lambda$

\item $|B(w_i,\delta)\cap\Gamma|\leqslant N$,
\end{itemize}
where $N$ is a constant that can be bounded by the number of
uniformly discrete sets that form $\Gamma$. That is, we define $\widetilde{\Gamma}$ so that $w_i$ is the only point of $\widetilde{\Gamma}$ in $B(w_i,\delta)$, every
$z_j\in\Gamma$ is contained in some $B(w_i,\delta)$, and each one of these balls only
contains a finite bounded number of points of $\Gamma$. Taking this into account we can write
$\Gamma=\cup_{i\in\N}\{w_i^1,\dotsc,w_i^{N_i}\}$, so that:
\begin{itemize}
\item $w_i^1=w_i\in \widetilde{\Gamma}$

\item $w_i^k\notin\widetilde{\Gamma}$ if $k\neq 1 $

\item $w_i^k\in B(w_i,\delta)$

\item $N_i\leqslant N\;\forall i$

\item The sets $\{w_i^k\}_{k=1}^{N_i}$ are disjoint.
\end{itemize}

With this notation we compute:
\begin{align*}
\sum_{z_j\in\Gamma}|&F(z_j)|^2=\sum_{w_i\in\widetilde{\Gamma}}\sum_{k=1}^{N_i}|F(w_i^k)|^2\\
=&\sum_{w_i\in\widetilde{\Gamma}}\left[\sum_{k=1}^{N_i}\left(|F(w_i^k)|^2-|F(w_i^1)|^2\right)+N_i|F(w_i^1)|^2\right]\\
\leqslant&N\sum_{w_i\in\widetilde{\Gamma}}|F(w_i)|^2+\sum_{w_i\in\widetilde{\Gamma}}
\sum_{k=1}^{N_i}\left(|F(w_i^k)|^2-|F(w_i^1)|^2\right).
\end{align*}

If we define $w_i^k=w_i^1$ for $N_i<k\leqslant N$ we can write:
\begin{align*}
\sum_{w_i\in\widetilde{\Gamma}}\sum_{k=1}^{N_i}&\left(|F(w_i^k)|^2-|F(w_i^1)|^2\right)\\
=&\sum_{w_i\in\widetilde{\Gamma}}\sum_{k=1}^N\left(|F(w_i^k)|^2-|F(w_i^1)|^2\right)\\
=&\sum_{k=1}^N\Biggl[\biggl(\sum_{w_i\in\widetilde{\Gamma}}|F(w_i^k)|^2\biggr)^{1/2}+
\biggl(\sum_{w_i\in\widetilde{\Gamma}}|F(w_i^1)|^2\biggr)^{1/2}\Biggr]\\
&\quad
\Biggl[\biggl(\sum_{w_i\in\widetilde{\Gamma}}|F(w_i^k)|^2\biggr)^{1/2}-
\biggl(\sum_{w_i\in\widetilde{\Gamma}}|F(w_i^1)|^2\biggr)^{1/2}\Biggr].
\end{align*}

We can be bound the first bracket by $2 B^{1/2}\|F\|$, since they
are partial sums of $\Gamma$. To bound the second bracket in
absolute value we apply \ref{gestab_sampling1} to the sets
$\{w_i^k\}_{i\in\N}$ and $\{w_i^1\}_{i\in\N}$ to get for each $k$
\begin{equation*}
\left|\left(\sum_{i\in\N}|F(w_i^k)|^2\right)^{1/2}-\left(\sum_{i\in\N}|F(w_i^1)|^2\right)^{1/2}\right|\leqslant
d_1^kd_2^k\|F\|,
\end{equation*}
where $d_1^k,\,d_2^k$ are defined by:
\begin{itemize}
\item $(d_1^k)^2=\sup_{i\in\N}\int_{\C}|k_{w_i^k-w_i^1}(z)-k(z)|\,d\mu(z)$

\item $(d_2^k)^2=\sup_{z\in\C}\sum_{i\in\N}|k_{w_i^k}(z)-k_{w_i^1}(z)|$.
\end{itemize}

We define now $d_1=\sup_k d_1^k$, $d_2=\sup_k d_2^k$ and bound independently of $k$:
\begin{equation*}
\left|\left(\sum_{i\in\N}|F(w_i^k)|^2\right)^{1/2}-\left(\sum_{i\in\N}|F(w_i^1)|^2\right)^{1/2}\right|
\leqslant d_1d_2\|F\|.
\end{equation*}

Therefore we obtain:
\begin{align*}
A\|F\|^2\leqslant&\sum_{j\in\N}|F(z_j)|^2\\
\leqslant&
N\sum_{i\in\N}|F(w_i)|^2+\left|\sum_{i\in\N}\sum_{k=1}^N\left(|F(w_i^k)|^2-|F(w_i^1)|^2\right)\right|\\
\leqslant& N\sum_{i\in\N}|F(w_i)|^2+2NB^{1/2}d_1d_2\|F\|^2.
\end{align*}

This implies
\begin{equation*}
\frac{A-2NB^{1/2}d_1d_2}{N}\|F\|^2\leqslant\sum_{i\in\N}|F(w_i)|^2.
\end{equation*}

As $N,\,A$ and $B$ are fixed, applying \ref{gcotasuperior} and
\ref{gdspetites} we see that there is $\delta$ small enough such
that
\begin{equation*}
A'\|F\|^2\leqslant\sum_{i\in\N}|F(w_i)|^2\leqslant B'\|F\|^2,
\end{equation*}
that is, $\widetilde{\Gamma}$ is a
sampling set.
\end{proof}

No we prove the existence of sampling sets, as we can find in \cite{FS2}.

\begin{lem}\label{gcompardiscr-con}
Let $\Gamma=\{z_j\}_{j\in\N}$ be such that for all $j$ there
exists an open set $V_j\subseteq\C$ so that $V_j\cap
V_k=\emptyset,\, j\neq k$ and $\C=\cup_{j\in\N}\overline{V_j}$,
with $z_j\in V_j$ and $\cup_{j\in\N} (V_j-z_j)\subseteq V$, with
$V$ compact and symmetrical. Then
\begin{equation*}
\Biggl|\|F\|-\biggl(\sum_{j\in\N}c_j|F(z_j)|^2\biggr)^{1/2}\Biggr|\leqslant
\widetilde{d}_1\widetilde{d}_2\|F\|,
\end{equation*}
where
\begin{itemize}
\item $\widetilde{d}_1^2=\sup_{z\in V}\int_{\C}
|k_z(w)-k(w)|\,dm(w)$

\item $\widetilde{d}_2^2=\sup_{w\in\C}\sum_{j\in\N}\int_{V_j}
|k_w(z)-k_w(z_j)| \,dm(z)$,
\end{itemize}
and $c_j$ is the area of $V_j$.
\end{lem}

\begin{proof}
Calculating directly as before, that is, looking at $\|F\|$ as the $l^2$-norm of the sequence $\bigl\{\bigl(\int_{V_j}|F(z)|^2\,dm(z)\bigr)^{\frac{1}{2}}\bigr\}$,
\begin{small}
\begin{equation*}
\Biggl|\|F\|-\biggl(\sum_{j\in\N}c_j|F(z_j)|^2\biggr)^{\frac{1}{2}}\Biggr|^2
\leqslant\sum_{j\in\N}\int_{V_j}\bigl||F(z)|-|F(z_j)|\bigr|^2\,dm(z).
\end{equation*}
\end{small}

By the invariance by translations this is equal to
\begin{equation*}
\sum_{j\in\N}\int_{V_j-z_j}\bigl||F_{-z_j}(z)|-|F_{-z_j}(0)|\bigr|^2dm(z).
\end{equation*}

We introduce the reproduction formula and use the Schwarz inequality to bound this by
\begin{align*}
\sum_{j\in\N}\int_{V_j-z_j}&\left[\int_{\C}|F_{-z_j}(w)|^2|k_z(w)-k(w)|dm(w)\right]\\
\cdot&\left[\int_{\C}|k_z(w)-k(w)|dm(w)\right]dm(z).
\end{align*}

We obtain the desired result bounding separately each part.
\end{proof}

\begin{thm}\label{gproudensessampling}
Let $H$ be the model space of a Gabor atom $g\in\mathcal{A}$.
There is $\delta$ such that if $\Gamma$ is a uniformly discrete set such
that $B(z,\delta)\cap\Gamma\neq\emptyset\;\forall z$ then $\Gamma$
is a sampling set for $H$.
\end{thm}

\begin{proof}
First we prove that there exists $\delta$ such that if $\Gamma=\{z_j\}_{j\in\N}$ is a uniformly discrete
set fulfilling the conditions of \ref{gcompardiscr-con} with $V$
contained in $B(0,\delta)$, and that $B(0,\alpha/2)\subset V_j-z_j$ for some $\alpha>0$, then $\Gamma$ is a sampling set for
$H$.

To see this we will bound $d_2'$ and to make $d_1'$ as small as we
want in \ref{gcompardiscr-con}. For $d_2'$, we fix $w\in\C$ and we
have that:
\begin{align*}
\sum_{j\in\N}\int_{V_j}|k_w(z)&-k_w(z_j)|\,dm(z)\\
 \leqslant&\sum_{j\in\N}\int_{V_j}|k_w(z)|\,dm(z)+\sum_{j\in\N}|V_j||k_w(z_j)|\\
\leqslant&\|k\|_1+\|Mk\|_1
\end{align*}
if $\delta\leqslant 1 $. Therefore
$d_2'\leqslant\bigl(\|k\|_1+\|Mk\|_1\bigr)^{1/2}$. To bound $d_1'$ we fix $z\in V$. As $k\in L^1(\C)$, there exists $R$ such that for every $z\in V$
\begin{equation*}
\int_{\C\setminus
B(0,R)}|k(w-z)|\,dm(w)\leqslant\frac{\varepsilon^2}{4}.
\end{equation*}

Now we write
\begin{align*}
\int_{\C} |&k_z(w)-k(w)|\,dm(w)\\
\leqslant& \int_{B(0,R)}\left|e^{2\pi ix(b-y)}k(w-z)- k(w)\right|\,dm(w)+\frac{\varepsilon^2}{2}\\
\leqslant& \int_{B(0,R)}\left|e^{2\pi ix(b-y)}-1\right||k(w-z)|\,dm(w)\\
&+\int_{B(0,R)}|k(w-z)-k(w)|\,dm(w)+\frac{\varepsilon^2}{2}.
\end{align*}

Now, as $w\in B(0,R)$ we have that $|b|<R$, and in an equivalent
way $|x|,|y|<\delta$. Therefore, if we take $\delta$ small enough,
we have that $\left|e^{2\pi
ix(b-y)}-1\right|<\frac{\|k\|_1\varepsilon^2}{4}$, which bounds
the first integral. We can make the second integral smaller than
$\frac{\varepsilon^2}{4}$ choosing $\delta$ small by the
continuity of the translation operator in $L^1(\C)$. Then we deduce that if $\delta$ is small enough we can
achieve $d_1'<\varepsilon$ for any $\varepsilon>0$. Therefore,
applying \ref{gcompardiscr-con} we see that:
\begin{equation*}
\Biggl|\|F\|-\biggl(\sum_{j\in\N}c_j|F(z_j)|^2\biggr)^{1/2}\Biggr|<\varepsilon(\|k\|_1+\|Mk\|_1)^{1/2}\|F\|
\end{equation*}
for small enough $\delta$. Choosing $\varepsilon$ so that
$\varepsilon\bigl(\|k\|_1+\|Mk\|_1\bigr)^{1/2}<1$ and noticing
that the conditions imply that the constants $c_j$, the area of
$V_j$, are bounded above and below because $B(0,\alpha/2)\subseteq
V_j-z_j \subseteq B(0,\delta)$, the proof of the first part is
finished.

Now we prove that we are in these conditions. We assume that
$\Gamma$ is an uniformly discrete set such that
$B(z,\delta)\cap\Gamma\neq\emptyset\forall z$ with the $\delta$
found before. Let $\alpha$ be the separation constant of $\Gamma$.
We want to see that $\exists V_j$ open, with $V_j\cap
V_k=\emptyset$ if $j\neq k$, $\C=\cup_{j\in\N}\overline{V_j}$,
$(V_j-z_j)\subseteq V$ compact and symmetric, and
$B(0,\alpha/2)\subseteq\cap_{j\in\N}(V_j-z_j)$. We define
$B_j=B(z_j,\delta)$ (from the definitions of $\alpha$ and $\delta$
it is clear that $2\delta>\alpha$). We now define the $V_j$:
\begin{itemize}
\item $V_1=B_1\setminus\cup_{j\neq1}b_j$

\item $V_k=B_i\setminus\left(\cup_{j\neq k}b_j\right)\bigcup\overline{
\left(\cup_{j=1}^{k-1}B_i\right)}$.
\end{itemize}

We observe that in each compact of $\C$ all these unions and
intersections are finite, and therefore the union is closed. This
says that the $V_j$ are all open sets, and for construction they
are disjointed. As $B_j\subseteq\cup_{k=1}^j\overline{V_k}$, and
every $z$ is in some $B_j$, we already have that
$\cup_{j\in\N}\overline{V_j}=\C$. As $V_j\subseteq B_j\Rightarrow
(V_j-z_j)\subseteq
(B_j-z_j)=B(0,\delta)\subseteq\overline{B(0,\delta)}$ compact and
symmetric, and also by construction, $B(0,\alpha/2)\subseteq
(V_j-z_j)$ and we obtain the last condition.
\end{proof}

\section{Sampling results for the wavelet transform.}

In this section we will give the same results that in the previous
section but for the wavelet transform. As mentioned in the
introduction, some results are already known after \cite{OS1},
although here we improve some of them and give explicit bounds.
Technically we need considering wavelets in $ \mathcal{MB}$, a
condition that in our understanding is required and missing in
\cite{OS1}. This set of functions plays the role of the
Feichtinger Algebra, but here, as opposed to the Gabor case, the
integrability of the kernel does not imply the integrability of
its maximal function. The proofs are very similar to those of the
former section, and we will point out only the substancial
diferences.

We establish first some notation. $H$ will always be the model
space of an admissible wavelet $\psi$ normalized so that
$\|\psi\|=1$ and
$\int_0^{\infty}\frac{|\widehat{\psi}(\xi)|^2}{\xi}\,d\xi=1$. We
remember that the balls and the separation conditions are with
respect the hyperbolic distance.

\begin{defi} A set $\Sigma=\{\sigma_j\}_{j\in\N}\subset \R\times\R^+$ is
said to be \emph{uniformly discrete} if there exists $\varepsilon$
such that $d(z_i,\sigma_j)>\varepsilon$ for $i\neq j$. The
supremum of such $\varepsilon$ is called the \emph{separation
constant} of $\Sigma$.
\end{defi}

\begin{prop}\label{unioseparades} If
$\Sigma$ is a sampling set for $H$ then $\Sigma$ is a finite union
of uniformly discrete sets.
\end{prop}

\begin{proof}
Same as \ref{gunioseparades}.
\end{proof}

\begin{lem}\label{acotdiscretKenL1}
Let $k$ be such that $Mk\in L^1(\R\times\R^+)$ and
$\Sigma\subset\C$ be a uniformly discrete set with separation
constant $\varepsilon$. Then
\begin{equation*}
\sum_{\sigma\in\Sigma}|k(z)|<\frac{(\sinh\frac{\varepsilon}{4})^{-2}}{4\pi}\|Mk\|_1.
\end{equation*}
\end{lem}

\begin{proof}
The same as in \ref{gacotdiscretKenL1}, but using the hyperbolic mesure. We suppose without losing generality that
$\frac{\varepsilon}{2}<1$. Then
\begin{equation*}
\sum_{\sigma\in\Sigma}|k(\sigma)|\leqslant\sum_{\sigma\in\Sigma}\frac{1}{|B(\sigma,\frac{\varepsilon}{2})|}
\int_{B(\sigma,\frac{\varepsilon}{2})}Mk(z)\,d\mu(z).
\end{equation*}
Since by hypothesis those balls are disjoint the lemma follows.
\end{proof}
\begin{obs}
Again notice that the bound only depends on the separation
constant of $\Sigma$ and that we are using that
$\varphi\in\mathcal{MB}$.
\end{obs}

\begin{prop}\label{cotasuperior}
Given a wavelet $\psi\in\mathcal{MB}$ and $\Sigma$ a uniformly
discrete set, there exists $B>0$ such that
\begin{equation*}
\sum_{\sigma\in\Sigma}|F(\sigma)|^2\leqslant B\|F\|^2\quad \forall
F\in H.
\end{equation*}
\end{prop}

\begin{proof}
The same as in \ref{gcotasuperior}.
\end{proof}

\begin{thm}
Let $H$ be the model space of an admissible wavelet
$\psi\in\mathcal{MB}$. Given $\Sigma=\{\sigma_j\}_{j\in\N}$ a
sampling set for $H$ there exists $\delta>0$ such that if
$\Gamma=\{w_j\}_{j\in\N}$ satisfies
$d(\sigma_j,w_j)<\delta\;\forall j$, $\Gamma$ is also a sampling
set.
\end{thm}

This result is deduced in a trivial way from the following pair of
lemmas.

\begin{lem}\label{estab_sampling1}
Let $\Sigma=\{\sigma_j\}_{j\in\N}$ and $\Gamma=\{w_j\}_{j\in\N}$
be two sets in $\R\times\R^+$, then:
\begin{equation*}
\left|\biggl(\sum_{j\in\N}|F(\sigma_j)|^2\biggr)^{\frac{1}{2}}-\biggl(\sum_{j\in\N}|F(w_j)|^2\biggr)^{\frac{1}{2}}\right|
\leqslant d_1d_2\|F\|\forall F\in H
\end{equation*}
where $d_1$ and $d_2$ are defined as:
\begin{itemize}
\item $d_1^2=\sup_j\int_{\R\times\R^+}|k(z^{-1}\cdot
\sigma_j)-k(z^{-1}\cdot w_j)|\,d\mu(z)$

\item $d_2^2=\sup_z\sum_{j\in\N}|k(z^{-1}\cdot \sigma_j)-k(z^{-1}\cdot
w_j)|$
\end{itemize}
\end{lem}

\begin{proof}
The same as in \ref{gestab_sampling1}.
\end{proof}

\begin{lem}\label{dspetites}
Let $\Sigma=\{\sigma_j\}_{j\in\N}$ be a uniformly discrete set, and assume
that $\psi\in\mathcal{MB}$. Then for every $\varepsilon$ there
exists $\delta$ such that if $\Gamma=\{w_j\}_{j\in\N}$ satisfies
$d(\sigma_j,w_j)\leqslant\delta\;\forall j$ then $d_1d_2\leqslant
\varepsilon$, with $d_1$ and $d_2$ defined as in
\ref{estab_sampling1}.
\end{lem}

\begin{proof}
First, we will see that we can make $d_1$ as small as we want if
$\Gamma$ is close enough to $\Sigma$. But in this case this is immediate, since
we only have to observe that if $|k(z)|$ is integrable
$|k(z^{-1})|$ also is, and that the translation operator is
uniformly continuous in $L^1(\R\times\R^+)$.

The bound $d_2$ and the rest of the proof is the same that in \ref{gdspetites}.
\end{proof}

\begin{thm}\label{samplingseparada}
Let $\Sigma=\{\sigma_j\}_{j\in\N}$ be a sampling set for the model
space of a wavelet $\psi\in\mathcal{MB}$. Then $\Sigma$ contains a
subset $\widetilde{\Sigma}\subseteq\Sigma$ such that
$\widetilde{\Sigma}$ is a sampling and uniformly discrete set.
\end{thm}

\begin{proof}
The proof is the same as in \ref{gsamplingseparada}, but using hyperbolic distance and balls.
\end{proof}

Now we prove the existence of sampling sets, recovering a result
in \cite{OS1}.

\begin{lem}\label{compardiscr-con}
Let $\Sigma=\{\sigma_j\}_{j\in\N}$ be such that for all $j$ there
exists an open set $V_j\subseteq\R\times\R^+$ so that $V_j\cap
V_k=\emptyset, j\neq k$ and
$\R\times\R^+=\cup_{j\in\N}\overline{V_j}$, with $\sigma_j\in V_j$
and $\cup_{j\in\N}\sigma_j^{-1}\cdot V_j\subseteq V$, with $V$
compact and symmetrical (around $i$). Then
\begin{equation*}
\Biggl|\|F\|-\biggl(\sum_{j\in\N}c_j|F(\sigma_j)|^2\biggr)^{\frac{1}{2}}\Biggr|\leqslant
\widetilde{d}_1\widetilde{d}_2\|F\|
\end{equation*}
where
\begin{itemize}
\item $\widetilde{d}_1^2=\sup_{\xi\in
V}\int_{\R\times\R^+}|k(z\cdot\xi)-k(z)|\,d\mu(z)$

\item
$\widetilde{d}_2^2=\sup_{w\in\R\times\R^+}\sum_{j\in\N}\int_{V_j}|k(w\cdot
z)-k(w\cdot \sigma_j)|\,d\mu(z)$
\end{itemize}
and $c_j$ is the (hyperbolic) area of $V_j$.
\end{lem}

\begin{proof}
The proof is the same that in \ref{gcompardiscr-con} and it can be found in \cite{OS1}.
\end{proof}

\begin{thm}\label{proudensessampling}
Let $H$ be the model space of a wavelet $\psi\in\mathcal{MB}$.
There is $\delta$ such that if $\Sigma$ is an uniformly discrete set such
that $B(z,\delta)\cap\Sigma\neq\emptyset\;\forall z$ then $\Sigma$
is a sampling set for $H$.
\end{thm}

\begin{proof}
The proof of the first part is here a bit different. We prove that
there is a $\delta$ such that if $\Sigma=\{\sigma_j\}_{j\in\N}$ is
a discrete set fulfilling the conditions of \ref{compardiscr-con}
with $V$ contained in $B(i,\delta)$, and that
$B(i,\alpha/2)\subseteq \sigma_j^{-1}\cdot V_j$ for some
$\alpha>0$, then $\Sigma$ is a sampling set for $H$.

To see this we will bound $\widetilde{d}_2$ and make $\widetilde{d}_1$ as small as we want in
\ref{compardiscr-con}. For $\widetilde{d}_2$, we fix $w\in\R\times\R^+$ and
we have that:
\begin{align*}
\sum_{j\in\N}\int_{V_j}\bigl|k(w\cdot z&)-k(w\cdot
\sigma_j)\bigr|\,d\mu(z)\\
\leqslant& \sum_{j\in\N}\int_{V_j}|k(w\cdot z)|\,d\mu(z)
+\sum_{j\in\N}|V_j||k(w\cdot \sigma_j)|\\
\leqslant&\|k\|_1+\|Mk\|_1
\end{align*}
if $\delta\leqslant 1 $. Therefore
$\widetilde{d}_2\leqslant(\|k\|_1+\|Mk\|_1)^{1/2}$. We note again
that we use the condition $\psi\in\mathcal{MB}$. By the
continuity of the translation operator in $L^1(\R\times\R^+)$ we
see that $\widetilde{d}_1$ can be made arbitrarily small if
$\delta$ is small enough. Therefore applying \ref{compardiscr-con}
we have that:
\begin{equation*}
\Biggl|\|F\|-\biggl(\sum_{j\in\N}c_j|F(\sigma_j)|^2\biggr)^{\frac{1}{2}}\Biggr|<\varepsilon\bigl(\|k\|_1+\|Mk\|_1\bigr)^{\frac{1}{2}}\|F\|
\end{equation*}
for small enough $\delta$. Choosing $\varepsilon$ so that
$\varepsilon(\|k\|_1+\|Mk\|_1)^{1/2}<1$, and noticing that the
conditions imply that the constants $c_j$, the hyperbolic areas of
the $V_j$, are bounded above and below because $B(i,\alpha/2)\subseteq
\sigma_j^{-1}\cdot V_j\subseteq B(i,\delta)$, the proof of the first part is finished.

The second part of the proof is identical to \ref{gproudensessampling} using hyperbolic distance and balls and we omit here.
\end{proof}

\end{document}